\documentclass[psamsfonts]{amsart}

\usepackage{float}
\usepackage{amsmath}
\usepackage{amsthm}
\usepackage{amssymb}
\usepackage{amscd}
\usepackage{amsfonts}
\usepackage{amsbsy}
\usepackage[noadjust]{cite} 
\usepackage{epsfig,afterpage}
\usepackage{paralist}
\usepackage{psfrag}
\usepackage{mathrsfs}       
\usepackage{verbatim}

\theoremstyle{plain}
\newtheorem {theorem} {Theorem} [section]

\theoremstyle{definition}

\newtheorem {example}[theorem] {Example}

\begin{document}

\title[Boundary Value Problems]
{Symbolic iterative Solution to Boundary Value Problem of Partial Differential Equations}

\author[H.~Semiyari ]{Hamid Semiyari$^1$}

\address{\noindent $^1$ Mathematics Department, American University, Washington, DC 20016, USA}

\email{semiyari@american.edu}



\begin{abstract}
In this article we introduce a simple straightforward and powerful method involving symbolic manipulation, Picard iteration, and auxiliary variables for approximating solutions of partial differential boundary value problems. The method is easy to implement, computationally efficient, and it is highly accurate. The output of the method is a function that approximates the exact solution.
\end{abstract}

\maketitle

\renewcommand{\thefootnote}{}
\footnote{${}^1$ Mathematics Department, American University, Washington, DC 20016, USA}

\section{Introduction}\label{s:intro}

In this article we introduce an efficient method to find analytical solution for boundary value problem of partial differential equations. The method is simple straightforward yet powerful. The method is based on the works of Paker-Sochacki \cite{PS} and Semiyari-Shafer \cite{SS}.

\section{The Algorithm for first order}\label{s:alg}

Let us consider a first order partial differential equation

     \begin{equation}\label{basic.bvp}
\begin{aligned}
u_t &= G(t,x, u,u_x)\\
u(a,x) &= f(x)\\
\end{aligned}
\end{equation}
and 
     \begin{equation}\label{bvp}
\begin{aligned}
u(t,b) &=\alpha(t)\\
\end{aligned}
\end{equation}
where $u_x\;=\;\frac{\partial u}{\partial x}$ and $u_t\;=\;\frac{\partial u}{\partial t}.$

Equation \eqref{basic.bvp} is equivalent to 

    \begin{equation}\label{e.baic}
\begin{aligned}
u(t,x)   &= f(x) \,+\, \int_a^t\,G(s,x, u(s,x),u_x(s,x))\,ds. \\
\end{aligned}
\end{equation}

If $u(t,x)$ is a solution of \eqref{basic.bvp} then \eqref{bvp} can be written as

   \begin{equation}\label{bvp-1}
\begin{aligned}
u(t,b)  \,-\, \alpha(t) =& 0. \\
\end{aligned}
\end{equation}
Hence, $u(t,x)$ is equivalent to
   \begin{equation}\label{bvp-2}
\begin{aligned}
u(t,x) &= u(t,x)\,-\,(u(t,b)  \,-\, \alpha(t)). \\
\end{aligned}
\end{equation}

The key idea in our proposed method to solve \eqref{basic.bvp} is that in a picard iteration scheme we use \eqref{bvp-2} to update the approximation of $u(t,x)$ so that the condition \eqref{bvp} be satisfied in each iteration. Thus, the iterates are,

   \begin{equation}\label{it.1}
\begin{aligned}
u^{[0]}(t,x) &= f(x), \\
\end{aligned}
\end{equation}
and
   \begin{equation}\label{it.2}
\begin{aligned}
u^{[k+1]}(t,x) &= f(x)\,+\, \int_a^t\,G(s,x, u^{[0]},u^{[k]}_{x})\,ds \\
u^{[k+1]}(t,x)  &= u^{[k+1]}(t,x)\,-\,( u^{[k+1]}(t,b)\,-\,\alpha(t)) \\
\end{aligned}
\end{equation}
where $u^{[k]}_{x}\;=\;\frac{\partial u^{[k]}}{\partial x}.$\\

This gives us the approximation solution to \eqref{basic.bvp} as a function of $t$ and $x$.

\subsection{Test Example}
\begin{example}
 Consider
    \begin{equation}\label{ex1}
 \begin{aligned}
u_t &= -u_x+2+t+x\\
u(0,x) &= f(x)\\
u(t,0) &= \alpha(t)\\
    \end{aligned}
\end{equation}
Where $f(x)=1+x$ and $\alpha(t)=1+t$. \\
The exact solution is given by
\[
 \begin{aligned}
u(t,x) &= (1+t)(1+x)\\
    \end{aligned}
\]

The equation \eqref{ex1} can be set up as follow
    \begin{equation}\label{N.ex1}
 \begin{aligned}
u_t &= -u_x+2+t+x\\
u(0,x) &= f(x)\\
    \end{aligned}
\end{equation}
and
    \begin{equation}\label{b.ex1}
 \begin{aligned}
u(t,0) &= \alpha(t)\\
    \end{aligned}
\end{equation}

The initial value is $u^{[0]}(t,x)=1+x$.

Thus the iterates are
\begin{equation} \label{e:iteration1}
\begin{aligned}
u^{[k+1]}(t,x)   &= 1+x \, - \, \int_0^t\,(u^{[k]}_x\,+\, 2\,+\,s\,+\,x)\,ds. \\
u^{[k+1]}(t,x)  &= u^{[k+1]}(t,x)\,-\,( u^{[k+1]}(t,0)\,-\,(1\,+\,t)) \\
\end{aligned}
\end{equation}
We obtain the exact solution by our very first iteration,
\begin{equation} 
\begin{aligned}
u^{[1]}(t,x)   &= (1+x) + \, \int_0^t\,-1+2+s+x\,ds \\
&=(1+x)+ t+\frac{t^2}{2}+xt\\
u^{[1]}(t,x) &=  u^{[1]}(t,x)\,-\,(u^{[1]}(t,0)\,-\,\alpha(t)) \\
&= 1+x+ t+\frac{t^2}{2}+xt\,-\,(1+ t+\frac{t^2}{2}\,-\,(1+t) \\
&=1+x+t+xt\\
&=(1+x)(1+t)
\end{aligned}
\end{equation}
Hence the absolute and relative errors are zero.
\end{example}

\begin{example}
In this example we consider a non-linear first order PDE.
\begin{equation}\label{basic.bvp.2}
 \begin{aligned}
A_2\,u_t + A_1\,u^{m+1}_x &= A_3\,u^j\,+\,A_4\,\exp(A_5\,u)+E(t,x).
\end{aligned}
\end{equation}
where $u_x\;=\;\frac{\partial u}{\partial x}$ and $u_t\;=\;\frac{\partial u}{\partial t}$
$$E(t,x)=\exp(t+y)\{A_2+A_1b^{-1}(m+1)\exp[m(t+y)]\}\,-\,A_3\exp[j(t+y)]\,-\,A_4[\exp[A_5\exp(t+y)]$$
 $f(x)=\exp(\frac{x}{b})$ and $\alpha(t)=\exp(t)$.\\
  
The initial values
\begin{equation}\label{iv.bvp2}
 \begin{aligned}
u(0,x) &= f(x)
\end{aligned}
\end{equation}
and boundary values
\begin{equation}\label{bv.bvp2}
 \begin{aligned}
u(t,0)  &= \alpha(t).
\end{aligned}
\end{equation}

The exact solution is
 \[
 \begin{aligned}
u(t,x) &= \exp(t+y),\quad y=\frac{x}{b}.\\
    \end{aligned}
\]  

The problem \eqref{basic.bvp.2}-\eqref{bv.bvp2} is written as 
 \begin{equation}\label{n.basic.bvp2}
 \begin{aligned}
u_t  &=\frac{1}{A_2}\Big(-A_1\,u^{m+1}_x+ A_3\,u^j\,+\,A_4\,\exp(A_5\,u)+E(t,x)\Big).\\
u(0,x) &= f(x)
\end{aligned}
\end{equation}
and
\begin{equation}\label{n.bv.bvp2}
 \begin{aligned}
u(t,0)  &= \alpha(t).
\end{aligned}
\end{equation}

We define the axillary variables to be
\begin{equation} \label{e:aux}
\begin{aligned}
v   &= \exp(t) \\
T   &= \exp(-u) \\
P   &= \exp(t+y) \\
R   &= \exp(-P) \\
\end{aligned}
\end{equation}
For ease of calculation we define $A_6\;=\;-A_5$.
Hence we will have a system of first order ode
\begin{equation} \label{e:system}
\begin{aligned}
u'  &= \frac{1}{A_2}\{(-A_1(u^{m+1})_x +A_3 u^j + A_4T^{A_6} + E\}\\
v'   &= v \\
T'   &= -u'T \\
P'   &= P \\
R'   &= -PR \\
E'   &= P\{A_2+A_1b^{-1}(m+1)^2P^m\}\,-\,A_3jP^j\,+\,A_4A_6PR^{A_6}
\end{aligned}
\end{equation}
where ``prime`` notation represents derivative with respect to $t$.\\
 The initial values or the first approximation for each variables are
\begin{equation} \label{e:aux}
\begin{aligned}
u^{[0]}    &= f\\
v^{[0]}   &= 1 \\
T^{[0]}   &= \exp(-f) \\
P^{[0]}   &= \exp(y) \\
R^{[0]}   &= \exp(-\exp(y)) \\
E^{[0]} &=\exp(y)\{A_2+A_1b^{-1}(m+1)\exp[my]\}\,\\
&\quad-\,A_3\exp(jy)\,-\,A_4\exp[-A_6\exp(y)]
\end{aligned}
\end{equation}
Let us define $$E_0\;=\;E[0].$$
The iterations would be

\begin{equation} \label{e:iteration1.bc}
\begin{aligned}
u^{[k+1]}    &= f  + \frac{1}{A_2}\int_0^t\,-A_1([u^{[k]}]^{m+1})_x +A_3 [u^{[k]}]^j + A_4[T^{[k]}]^{A_6} + E^{[k]} \,ds \\
v^{[k+1]}   &= 1  + \int_0^t\,v^{[k]}\,ds \\
T^{[k+1]}   &= \exp(-f)  - \int_0^t\, u^{[k+1]}T^{[k]}\,ds \\
P^{[k+1]}   &= \exp(y)  + \int_0^t\, P^{[k]}\,ds \\
R^{[k+1]}   &= \exp(-\exp(y))  - \int_0^t\,P^{[k]}R^{[k]} \,ds \\
E^{[k+1]}  &=E_0  + \int_0^t\,P^{[k]}\{A_2+A_1b^{-1}(m+1)^2[P^{[k]}]^m\}\,-\,A_3k[P^{[k]}]^j\,+\,A_4A_6P^{[k]}[R^{[k]}]^{A_6}\,ds \\
u^{[k+1]}(t,x)   &= u^{[k+1]}(t,x) \,-\,(u^{[k+1]}(t,0)\,-\,\alpha(t)) \\
\end{aligned}
\end{equation}

Let us take the initial data suggested in \cite{Y}, where
$m=1,\;b=1,\;A_1=1,\;A_2=1$. We have considered the several cases. After 4 iterations the error plot and maximum error provided
\begin{enumerate}
\item  {\bf Case 1:} $j=1,\;A_3=0,\;A_4=0,\;A_5=0$.\\

Absolute and Relative plots of the error in the final approximation $u^{[5]}(t,x)$ are shown in Figure \ref{Case 1} and \ref{Case 1.1} respectively. 

The maximum error rounded to 5 decimal digit; after 4 iterations the maximum error is $\epsilon \approx 0.00439$

\begin{figure}[H]
\centering
\includegraphics[width=3in]{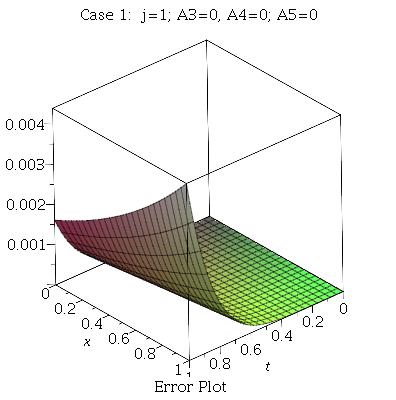}
\caption{Error plot, $|u(t,x)-u_{approx}|$, with $n=4$ iterations. $u(t,x)$ represents the exact solution and $u_{approx}$ represents the approximate solution.}\label{Case 1}
\end{figure}

\begin{figure}[H]
\centering
\includegraphics[width=3in]{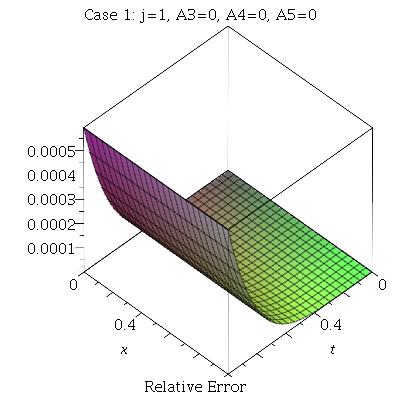}
\caption{Relative Error plot, $\frac{|u(t,x)-u_{approx}|}{|u(t,x)|}$, with $n=4$ iterations. $u(t,x)$ represents the exact solution and $u_{approx}$ represents the approximate solution.}\label{Case 1.1}
\end{figure}


\item  {\bf Case 2:} $j=1,\;A_3=-1,\;A_4=0,\;A_5=0$.\\

Absolute and Relative plots of the error in the final approximation $u^{[5]}(t,x)$ are shown in Figure \ref{Case 2} and \ref{Case 2.1} respectively. 

\begin{figure}[H]
\centering
\includegraphics[width=3in]{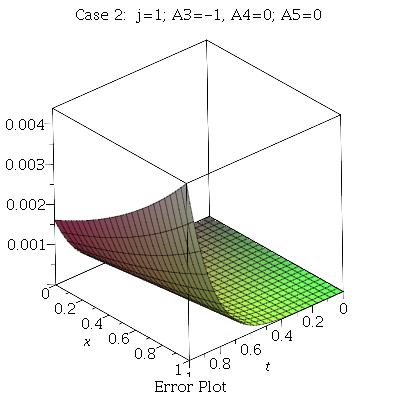}
\caption{Error plot, $|u(t,x)-u_{approx}|$, with $n=4$ iterations. $u(t,x)$ represents the exact solution and $u_{approx}]$ represents the approximate solution.}\label{Case 2}
\end{figure}

\begin{figure}[H]
\centering
\includegraphics[width=3in]{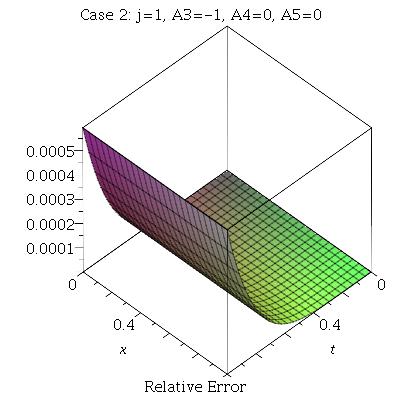}
\caption{Relative Error plot.}\label{Case 2.1}
\end{figure}


\item  {\bf Case 3:} $j=2,\;A_3=-1,\;A_4=0,\;A_5=0$.\\

Absolute and Relative plots of the error in the final approximation $u^{[5]}(t,x)$ are shown in Figure \ref{Case 3} and \ref{Case 3.1} respectively. 

\begin{figure}[H]
\centering
\includegraphics[width=3in]{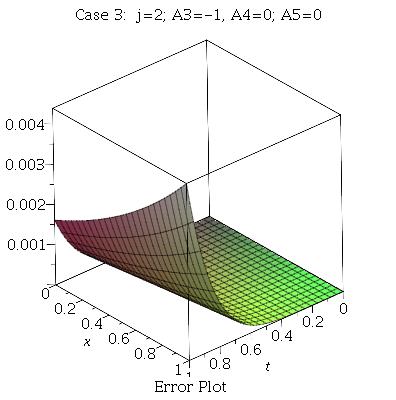}
\caption{Error plot, $|u(t,x)-u_{approx}|$, with $n=4$ iterations. $u(t,x)$ represents the exact solution and $u_{approx}$ represents the approximate solution.}\label{Case 3}
\end{figure}

\begin{figure}[H]
\centering
\includegraphics[width=3in]{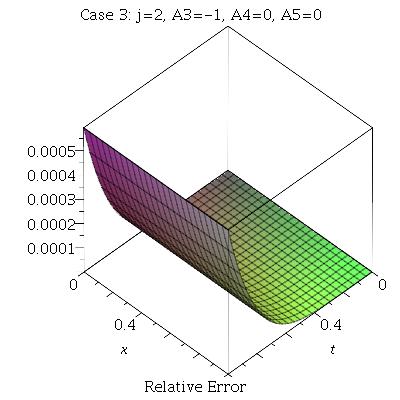}
\caption{Relative Error plot.}\label{Case 3.1}
\end{figure}


\item  {\bf Case 4:} $j=2,\;A_3=1,\;A_4=0,\;A_5=0$.\\

Absolute and Relative plots of the error in the final approximation $u^{[5]}(t,x)$ are shown in Figure \ref{Case 4} and \ref{Case 4.1} respectively. 

\begin{figure}[H]
\centering
\includegraphics[width=3in]{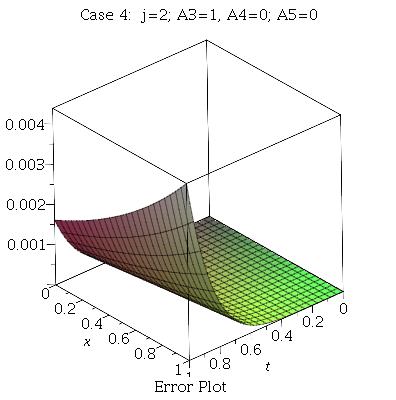}
\caption{Error plot, $|u(t,x)-u_{approx}|$, with $n=4$ iterations. $u(t,x)$ represents the exact solution and $u_{approx}$ represents the approximate solution.}\label{Case 4}
\end{figure}

\begin{figure}[H]
\centering
\includegraphics[width=3in]{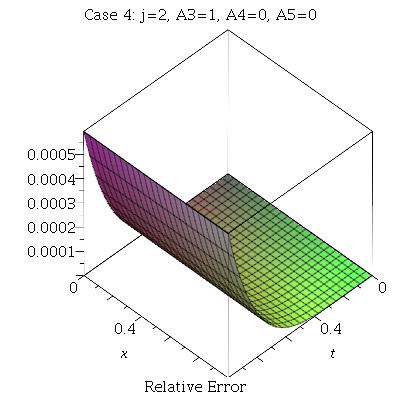}
\caption{Relative Error plot.}\label{Case 4.1}
\end{figure}


\item  {\bf Case 5:} $j=1,\;A_3=1,\;A_4=0,\;A_5=0$.\\

Absolute and Relative plots of the error in the final approximation $u^{[5]}(t,x)$ are shown in Figure \ref{Case 5} and \ref{Case 5.1} respectively. 

\begin{figure}[H]
\centering
\includegraphics[width=3in]{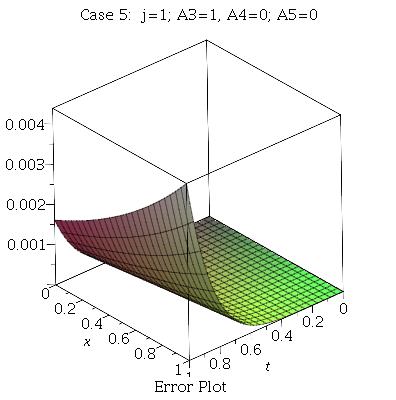}
\caption{Error plot, $|u(t,x)-u_{approx}|$, with $n=4$ iterations. $u(t,x)$ represents the exact solution and $u_{approx}$ represents the approximate solution.}\label{Case 5}
\end{figure}

\begin{figure}[H]
\centering
\includegraphics[width=3in]{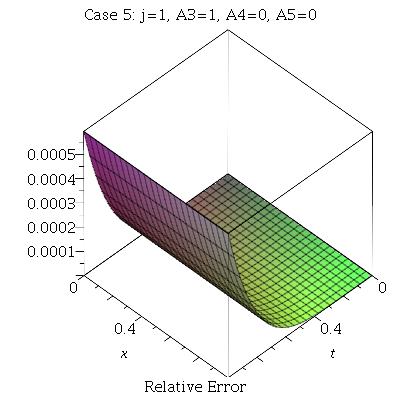}
\caption{Relative Error plot.}\label{Case 5.1}
\end{figure}


\item  {\bf Case 6:} $j=1,\;A_3=0,\;A_4=1,\;A_5=-1$.\\

Absolute and Relative plots of the error in the final approximation $u^{[5]}(t,x)$ are shown in Figure \ref{Case 6} and \ref{Case 6.1} respectively. 

\begin{figure}[H]
\centering
\includegraphics[width=3in]{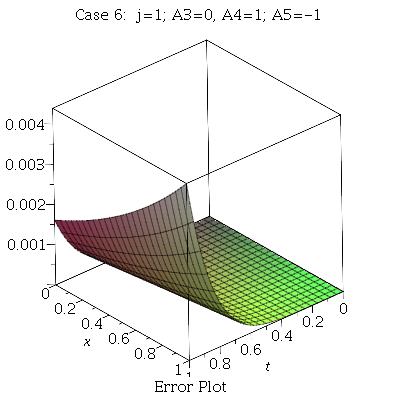}
\caption{Error plot, $|u(t,x)-u_{approx}|$, with $n=4$ iterations. $u(t,x)$ represents the exact solution and $u_{approx}$ represents the approximate solution.}\label{Case 6}
\end{figure}

\begin{figure}[H]
\centering
\includegraphics[width=3in]{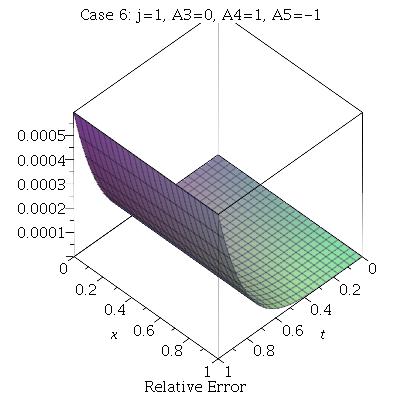}
\caption{Relative Error plot.}\label{Case 6.1}
\end{figure}

\item  {\bf Case 7:} $j=1,\;A_3=1,\;A_4=1,\;A_5=-1$.\\

Absolute and Relative plots of the error in the final approximation $u^{[5]}(t,x)$ are shown in Figure \ref{Case 7} and \ref{Case 7.1} respectively. 

\begin{figure}[H]
\centering
\includegraphics[width=3in]{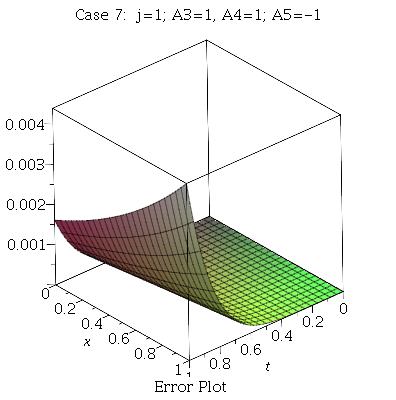}
\caption{Error plot, $|u(t,x)-u_{approx}|$, with $n=4$ iterations. $u(t,x)$ represents the exact solution and $u_{approx}$ represents the approximate solution.}\label{Case 7}
\end{figure}

\begin{figure}[H]
\centering
\includegraphics[width=3in]{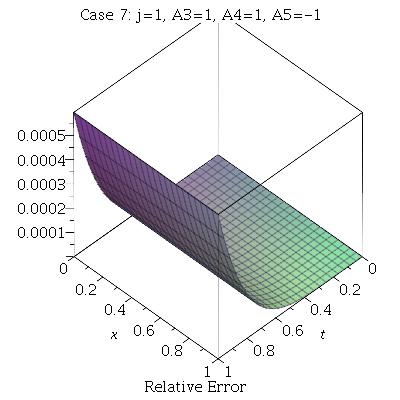}
\caption{Relative Error plot.}\label{Case 7.1}
\end{figure}

\end{enumerate}
\end{example}
\section{The Algorithm for second order}\label{s:alg2}

 Consider a second order PDE,

     \begin{equation}\label{basic.bvp1}
u_{xx}\;=\;G(t,x, u,u_x,u_x, u_{tt}, u_{tx}).
\end{equation}

Similar to the first order, we will set up a system of the first order Ordinary Differential Equations, ODE, by using auxiliary variables and then we will satisfy the boundary condition(s) iteratively.

\subsection{Test Examples}
\begin{example}
 Consider
    \begin{equation}\label{s.ex.1}
    \begin{aligned}
u_{tt}&=-u_{xx}
\end{aligned}
\end{equation}
Where $u(t,a)=\exp(t)$, $u(t,b)=0$, $u(a,x)=\cos x$ and $u_t(a,x)=\cos x$. \\
The exact solution is given by
\[ u(t,x)=\exp(t) \, \cos x\]
Let $\alpha=\exp(t)$, $\beta=0$, $a=0$, and $b=\frac{\pi}{2}$. 

The equation \eqref{s.ex.1} with its conditions can be written as
   \begin{equation}\label{e.ex.1}
    \begin{aligned}
u_{tt}&=-u_{xx}\\
u(a,x)&=\cos x\\
u_t(a,x)&=\cos x\\
\end{aligned}
\end{equation}
and 
  \begin{equation}\label{c.ex.1}
    \begin{aligned}
u(t,a)&=\exp(t)\;(=\alpha)\\
u(t,b)&=0\;(=\beta)\\
\end{aligned}
\end{equation}
If $u(t,x)$ is a solution to \eqref{e.ex.1} then by \eqref{c.ex.1} we can represent $u(t,x)$ as 
  \begin{equation}\label{c.ex1}
    \begin{aligned}
u(t,x)&= u(t,x) \,-\,\frac{x-a}{b-a}(u(t,b)\,-\,\beta)-\,\frac{b-x}{b-a}(u(t,a)\,-\,\alpha) \\
\end{aligned}
\end{equation}

We define the axillary variables to be
\begin{equation} \label{s:aux1}
\begin{aligned}
v   &= u_t \\
U  &= \exp(t).
\end{aligned}
\end{equation}
Hence we will have a system of first order ode
\begin{equation} \label{s:system1}
\begin{aligned}
u'  &= v\\
v'   &= -u_{xx} \\
U'  &= U\\
\end{aligned}
\end{equation}
where ``prime`` notation means derivative with respect to variable $t$.

The initial conditions are

 \begin{equation} \label{s:int1}
\begin{aligned}
u^{[0]}    &= \cos x\; (=u_0)\\
v^{[0]}   &= \cos x \;(=v_0) \\
U^{[0]}   &= 1\; (=U_0) \\
\end{aligned}
\end{equation}

 To solve \eqref{e.ex.1} is that in a picard iteration scheme we use \eqref{c.ex1} to update the approximation of $u(t,x)$ so that the condition \eqref{c.ex.1} be satisfied in each iteration. Thus, the iterates are

\begin{equation} \label{s:iteration1.bc.11}
\begin{aligned}
u^{[k+1]}(t,x)    &= u_0  + \int_0^t\,v^{[k]}(s,x) \,ds \\
v^{[k+1]}(t,x)   &= v_0  - \int_0^t\,u^{[k]}_{xx}(s,x)\,ds \\
U^{[k+1]}(x)   &= U_0+ \int_0^t\, U^{[k]}(s)\,ds \\
u^{[k+1]}(t,x)   &= u^{[k+1]}(t,x) \,-\,\frac{x-a}{b-a}(u^{[k+1]}(t,b)\,-\,\beta)-\,\frac{b-x}{b-a}(u^{[k+1]}(t,a)\,-\,\alpha). \\
\end{aligned}
\end{equation}

Absolute and Relative plots of the error in the final approximation $u^{[5]}(t,x)$ are shown in Figure \ref{p3} and \ref{p3-1} respectively. 

The maximum error rounded to 5 decimal digit; after 4 iterations the maximum error is $\epsilon \approx 0.0003$
\begin{figure}[H]
\centering
\includegraphics[width=3in]{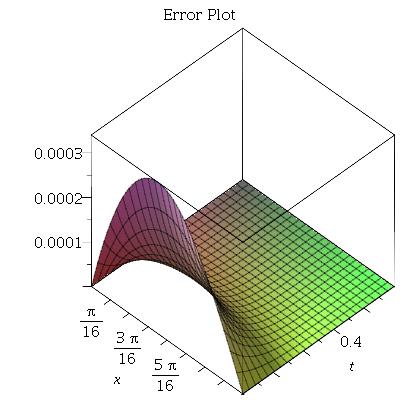}
\caption{Error plot, $|u(t,x)-u_{approx}|$, with $n=4$ iterations. $u(t,x)$ represents the exact solution and $u_{approx}$ represents the approximate solution.}\label{p3}
\end{figure}

\begin{figure}[H]
\centering
\includegraphics[width=3in]{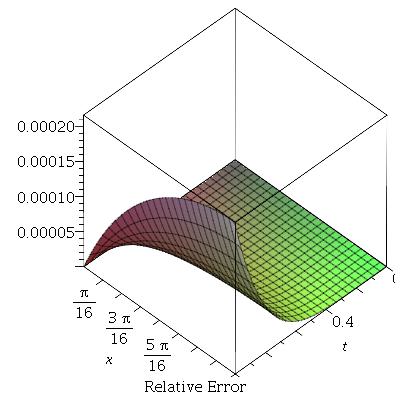}
\caption{Relative Error plot, $\frac{|u(t,x)-u_{approx}|}{|u(t,x)|}$, with $n=4$ iterations. $u(t,x)$ represents the exact solution and $u_{approx}$ represents the approximate solution.}\label{p3-1}
\end{figure}

\end{example}

\begin{example}{\bf Sine Gorden.}\\
 Consider
    \begin{equation}\label{s.ex2}
    \begin{aligned}
u_{xx}  &= u_{tt} - \sin u
\end{aligned}
\end{equation}
Where $u(t,a)=f(t)$, $u_x(t,a)=g(t)$, $u(a,x)=\alpha$ and $u(b,x)=\beta$. \\
The exact solution is given by \cite{Wolf}
\[ u(t,x)=-4\, \tan^{-1}\Big(\frac{m}{\sqrt{1-m^2}}\frac{\sin \sqrt{1-m^2}t}{\cosh (mx)}\Big)\]
Let $\alpha=0$, $\beta=-4\, \tan^{-1}\Big(\frac{m}{\sqrt{1-m^2}}\frac{\sin \sqrt{1-m^2}}{\cosh (mx)}\Big)$, $a=0$, and $b=1$. \[m^2\;<\;1\]

The equation \eqref{s.ex2} with its conditions can be written as
\begin{equation}\label{n.s.ex2}
    \begin{aligned}
u_{xx}  &= u_{tt} - \sin u\\
u(t,a)  &= f(t)\\
u_x(t,a)  &= g(t)
\end{aligned}
\end{equation}

and 

\begin{equation}\label{c.s.ex2}
    \begin{aligned}
u_{xx}  &= u_{tt} - \sin u\\
u(a,x)  &= \alpha(x)\\
u(b,x) &=  \beta(x)
\end{aligned}
\end{equation}

We define the axillary variables to be
\begin{equation} \label{s:aux2}
\begin{aligned}
v   &= u_x \\
U &= \sin u \\
V &= \cos u \\
\end{aligned}
\end{equation}
Hence we will have a system of first order ode
\begin{equation} \label{s:system2}
\begin{aligned}
u'  &= v\\
v'   &= u_{tt} - U \\
U'  &= vV\\
V'  &= -vU\\
\end{aligned}
\end{equation}
where the ``prime`` notation represent derivative with respect to $x$.\\
The initial values or the first approximation for each variables are
\begin{equation} \label{s:int2}
\begin{aligned}
u^{[0]}    &= f(t)\; (=u_0)\\
v^{[0]}   &= g(t) \;(=v_0) \\
U^{[0]}   &= \sin u_0\; (=U_0) \\
V^{[0]}   &= \cos V_0\; (=V_0) \\
\end{aligned}
\end{equation}

The iterations would be

\begin{equation} \label{s:iteration1.bc2}
\begin{aligned}
u^{[k+1]}(t,x)    &= u_0  + \int_0^x\,v^{[k]}(t,s) \,ds \\
v^{[k+1]}(t,x)   &= v_0  - \int_0^x\,u^{[k]}_{xx}(t,s)-\sin (u(t,s))\,ds \\
U^{[k+1]}(t,x)   &= U_0+ \int_0^x\, v^{[k]}V^{[k]}(s)\,ds \\
V^{[k+1]}(t,x)   &= V_0- \int_0^x\, v^{[k]}U^{[k]}(s)\,ds \\
u^{[k+1]}(t,x)   &= u^{[k+1]}(t,x) \,-\,\frac{t-a}{b-a}(u^{[k+1]}(b,x)\,-\,\beta)-\,\frac{b-t}{b-a}(u^{[k+1]}(a,x)\,-\,\alpha) \\
\end{aligned}
\end{equation}

Absolute and Relative plots of the error in the final approximation $u^{[5]}(t,x)$ are shown in Figure \ref{p4.1} and \ref{p4.1-1} respectively with $m=0.1$. The maximum error  approximated to $\epsilon \approx 0.010$ 

\begin{figure}[H]
\centering
\includegraphics[width=3in]{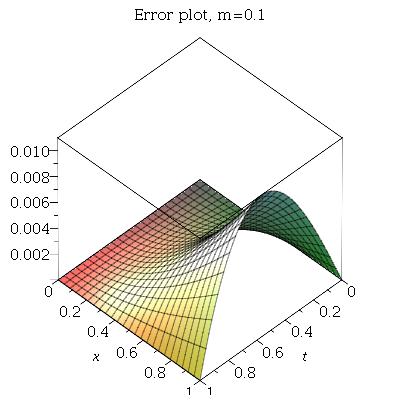}
\caption{Error plot, $|u(t,x)-u_{approx}|$, with $n=4$ iterations. $u(t,x)$ represents the exact solution and $u_{approx}$ represents the approximate solution.}\label{p4.1}
\end{figure}

\begin{figure}[H]
\centering
\includegraphics[width=3in]{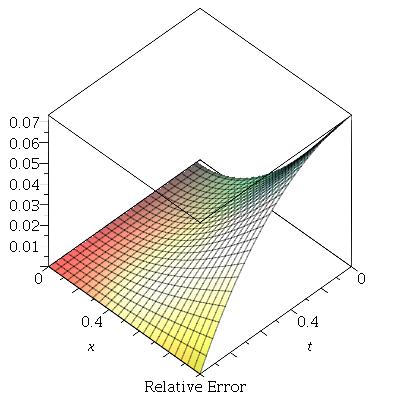}
\caption{Relative Error plot, $\frac{|u(t,x)-u_{approx}|}{|u(t,x)|}$, with $n=4$ iterations. $u(t,x)$ represents the exact solution and $u_{approx}$ represents the approximate solution.}\label{p4.1-1}
\end{figure}


Absolute and Relative plots of the error in the final approximation $u^{[5]}(t,x)$ are shown in Figure \ref{p4.2} and \ref{p4.2-1} respectively with $m=0.5$. The maximum error  approximated to $\epsilon \approx 0.05$ 

\begin{figure}[H]
\centering
\includegraphics[width=3in]{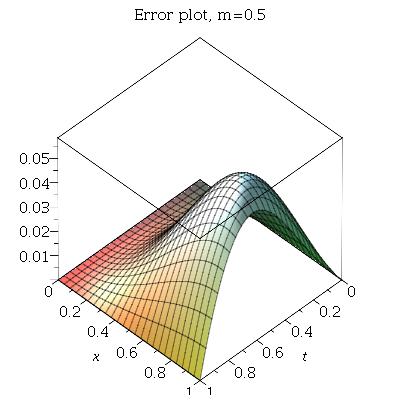}
\caption{Error plot, $|u(t,x)-u_{approx}|$, with $n=4$ iterations. $u(t,x)$ represents the exact solution and $u_{approx}$ represents the approximate solution.}\label{p4.2}
\end{figure}

\begin{figure}[H]
\centering
\includegraphics[width=3in]{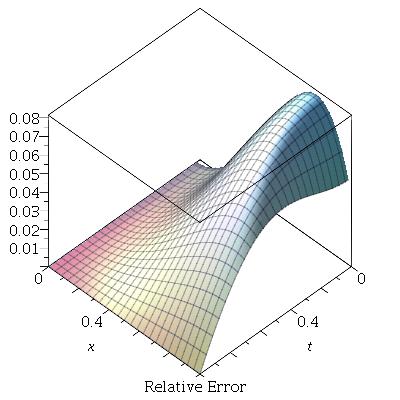}
\caption{Relative Error plot.}\label{p4.2-1}
\end{figure}


Absolute and Relative plots of the error in the final approximation $u^{[5]}(t,x)$ are shown in Figure \ref{p4.3} and \ref{p4.3-1} respectively with $m=0.9$. The maximum error  approximated to $\epsilon \approx 0.12$ 

\begin{figure}[H]
\centering
\includegraphics[width=3in]{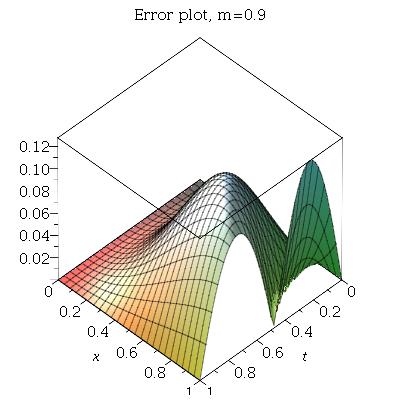}
\caption{Error plot, $|u(t,x)-u_{approx}|$, with $n=4$ iterations. $u(t,x)$ represents the exact solution and $u_{approx}$ represents the approximate solution.}\label{p4.3}
\end{figure}

\begin{figure}[H]
\centering
\includegraphics[width=3in]{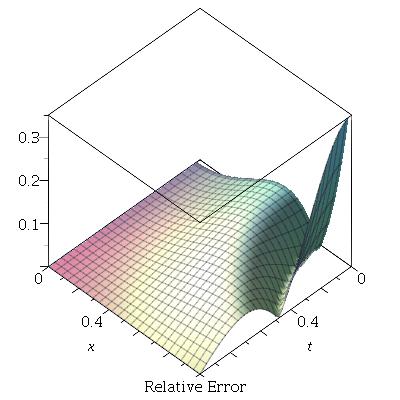}
\caption{Relative Error plot.}\label{p4.3-1}
\end{figure}

\end{example}


\begin{example}\label{ex5}

 Consider
    \begin{equation}\label{s.ex5}
    \begin{aligned}
u_{xx} &=  -2u_t\,u_x
\end{aligned}
\end{equation}
Where $u(t,a)=\alpha$ and $u(t,b)=\beta$ and $u(a,x)=\frac{2}{x+1}$. 
where $\alpha=t+2$, $\beta=\frac{2+t}{2}$, $a=0$ and $b=1$. \\
The exact solution is given by
\[
\begin{aligned}
u(t,x)   &= \frac{2+t}{1+x} \\
\end{aligned}
\]
The equation \eqref{s.ex5} with its condition is written as
    \begin{equation}\label{n.s.ex5}
    \begin{aligned}
u_{xx} &=  -2u_t\,u_x\\
u(t,0) &= \alpha(t)\\
u(t,1) &= \beta(t)\\
\end{aligned}
\end{equation}
and 
  \begin{equation}\label{c.s.ex5}
    \begin{aligned}
u(0,x) &= \frac{2}{x+1}\\
\end{aligned}
\end{equation}

We define the axillary variables to be
\begin{equation} \label{s:aux5}
\begin{aligned}
v   &= u_x \\
\end{aligned}
\end{equation}
Hence we will have a system of first order ode
\begin{equation} \label{s:system5}
\begin{aligned}
u'  &= v\\
v'   &= -2u_tv \\
\end{aligned}
\end{equation}
where the ``prime`` notation represents derivative with respect to $x$.\\

With initial values
\begin{equation} \label{s:initial5}
\begin{aligned}
u(t,0)  &= \alpha\\
v(t,0)  &= \gamma\\
\end{aligned}
\end{equation}
where $\gamma$ is unknown. However, it can be approximated by Semiyari-Shafer's method \cite{SS}.

 The first approximation for each variables are
\begin{equation} \label{s:aux}
\begin{aligned}
u^{[0]}    &= t+2\\
v^{[0]}   &= \frac{\beta-\alpha}{b-a} \\
\gamma^{[0]} &=\frac{1}{b-a}(\beta-\alpha + 2\int_0^1\,(b-s) u^{[0]}_t(t,s)\,v^{[0]}(t,s)\,ds)
\end{aligned}
\end{equation}

The iterations would be

\begin{equation} \label{s:iteration1.bc}
\begin{aligned}
u^{[k+1]}(t,x)    &= \alpha + \int_0^x\,v^{[k]}(t,s) \,ds \\
v^{[k+1]}(t,x)   &= \gamma^{[k]}  - 2\int_0^x\,l u^{[k]}_{t}(t,s)\,v^{[k]}(t,s)\,ds \\
\gamma^{[k+1]}   &= \frac{1}{b-a}(\beta-\alpha + 2\int_0^1\,(b-s) u^{[k]}_{t}(t,s)\,v^{[k+1]}(t,s)\,ds) \\
u^{[k+1]}(t,x)    &= u^{[k+1]}(t,x)    - (u^{[k+1]}(0,x) -\frac{2}{x+1}).\\    
\end{aligned}
\end{equation}
 After 4 iterations the exact $\gamma=-(t+2)$ approximated with 5 decimal digits to $\gamma^{[5]}=-(0.99673\,t+1.99346)$ 
 and the maximum error approximated to $\epsilon \approx 0.010$
The error in the final approximation $u^{[5]}(t,x)$ is shown in Figure \ref{p5}

\begin{figure}[H]
\centering
\includegraphics[width=3in]{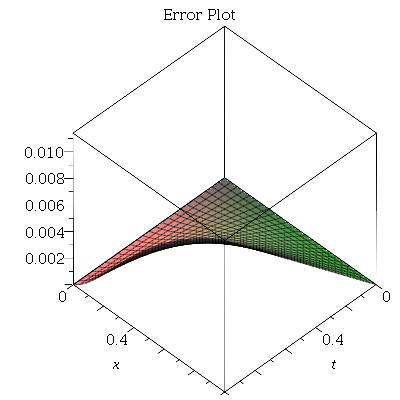}
\caption{Error plot, $|u(t,x)-u_{approx}|$, with $n=4$ iterations. $u(t,x)$ represents the exact solution and $u_{approx}$ represents the approximate solution.}\label{p5}
\end{figure}

The relative error is shown in Figure \ref{p5-1}

\begin{figure}[H]
\centering
\includegraphics[width=3in]{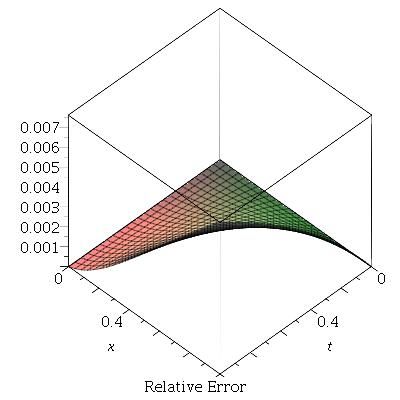}
\caption{The Relative Error plot, $\frac{|u(t,x)-u_{approx}|}{|u(t,x)|}$, with $n=4$ iterations. $u(t,x)$ represents the exact solution and $u_{approx}$ represents the approximate solution.}\label{p5-1}
\end{figure}

\end{example}

\begin{example} Let us revisit example \eqref{ex5}. The problem can be rewritten as

    \begin{equation}\label{s.ex6}
\begin{aligned}
u_t  &= -\frac{u_{xx}}{2u_x}
\end{aligned}
\end{equation}
Where $u(t,a)=\alpha$ and $u(t,b)=\beta$ and $u(0,x)=f$. 
where $\alpha=t+2$, $\beta=\frac{2+t}{2}$, $f=\frac{2}{x+1}$, $a=0$ and $b=1$. \\
The exact solution is given by
\[
\begin{aligned}
u(t,x)  &= \frac{2+t}{1+x}
\end{aligned}
\]

 The problem \eqref{s.ex6} with its conditions is written as
     \begin{equation}\label{n.s.ex6}
\begin{aligned}
u_t  &= -\frac{u_{xx}}{2u_x}\\
u(0,x) &= f
\end{aligned}
\end{equation}
and 
    \begin{equation}\label{c.s.ex6}
\begin{aligned}
u(t,0)  &= \alpha(t)\\
u(t,1) &= f\beta(t)
\end{aligned}
\end{equation}
if $u(t,x)$ is a solution of \eqref{n.s.ex6} then it is equivalent to 

\begin{equation}
\begin{aligned}
u(t,x)    &= u(t,x)    - \frac{x-a}{b-a}(u(t,b) -\beta)- \frac{b-x}{b-a}(u(t,a) -\alpha).\\      
\end{aligned}
\end{equation}

The system of first order ode
\begin{equation} 
\begin{aligned}
u_t  &= -\frac{u_{xx}}{2u_x}
\end{aligned}
\end{equation}
With initial values
\begin{equation} 
\begin{aligned}
u(0,x)  &= f\\
\end{aligned}
\end{equation}

 The first approximation is
\begin{equation} 
\begin{aligned}
u^{[0]}    &= f\\
\end{aligned}
\end{equation}

The iterations would be

\begin{equation}
\begin{aligned}
u^{[k+1]}(t,x)    &= f - \int_0^t\,\frac{u^{[k]}_{xx}}{2u^{[k]}_x}(t,s) \,ds \\
u^{[k+1]}(t,x)    &= u^{[k+1]}(t,x)    - \frac{x-a}{b-a}(u^{[k+1]}(t,b) -\beta)- \frac{b-x}{b-a}(u^{[k+1]}(t,a) -\alpha).\\      
\end{aligned}
\end{equation}
After 2 iterations the $u^{[3]}(t,x)=\frac{t+2}{x+1}$ which is same as the exact solution. Hence the Absolute and Relative errors are zero.

\end{example}

\section{Conclusion}\label{s:concl}

We have introduced a new analytic method for solving boundary value problem of partial differential equations. The method is easy to implement, computationally efficient, and it is highly accurate. The output of my method is a function that approximates the exact solution, whereas in the current methods the output is a sequence of n+1 points that approximate the values of the unknown solution at n+1 t-values and we are able to compute the error only at these n+1 points and then must interpolate between them.

\section{Acknowledgement}\label{s:Ack}
We thank James Sochacki, professor of Mathematics at James Madison University, for sharing his pearls of wisdom with us during the course of this research that greatly improved the manuscript.

\end{document}